\documentstyle[12pt]{article}
  \newcommand{\const}{\rm const}
  
  \newcommand{\supp}{\rm supp}

  \newcommand{\Dom}{\rm  Dom}

\begin{document}

   \begin{center}

  \   {\bf  A counterexample for equivalence result between tails } \par

\vspace{4mm}

 {\bf behavior and Grand Lebesgue Spaces norms }\\

\vspace{4mm}

  \   {\bf Ostrovsky E., Sirota L.}\\

\vspace{4mm}

 Israel,  Bar-Ilan University, department  of Mathematic and Statistics, 59200, \\

\vspace{4mm}

e-mails:\ eugostrovsky@list.ru \\
sirota3@bezeqint.net \\

\vspace{5mm}

  {\bf Abstract} \\

\vspace{4mm}

 \end{center}

 \  We bring in this short report  (counter -) examples on order to show a difference between tails behavior and Grand Lebesgue Spaces norms
 for some classes of random variables. \par

\vspace{4mm}

 \ {\it Key words and phrases:} Random variable and random vector (r.v.), centered (mean zero) r.v., saddle-point method, examples and counterexamples,
tail and  bilateral tail estimates, rearrangement invariant Banach space of  random variables, tail of distribution, Lorentz norms and spaces,
moments, Lebesgue-Riesz, Orlicz and Grand Lebesgue Spaces (GLS);  slowly varying functions, Tchebychev-Markov inequality, tail function,
Young-Fenchel transform, theorem and inequality of Fenchel-Moreau, Young-Orlicz
function, norm, Markov-Tchernov's estimate, non-asymptotical estimates, Cramer's condition. \par

\vspace{4mm}


\vspace{5mm}

\section{ Definitions.  Notations. Previous results.  Statement of problem.}

 \vspace{4mm}

 \begin{center}

\ {\bf A. \ $ \  B(\phi)  \  $ spaces.} \par

 \end{center}

 \ Let $ \  ( \ \Omega = \{ \omega\}, \ F, {\bf P} \ )  \ $ be certain sufficiently rich probability space. \par

 \ Let also $ \phi = \phi(\lambda), \lambda \in (-\lambda_0, \lambda_0), \ \lambda_0 =
\const \in (0, \infty] $ be certain even strong convex which takes positive values for positive arguments twice continuous
differentiable function, briefly: Young-Orlicz function, such that
$$
 \phi(0) = \phi'(0) = 0, \ \phi^{''}(0) > 0, \ \lim_{\lambda \to \lambda_0} \phi(\lambda)/\lambda = \infty. \eqno(1.1)
$$

 \ For instance: $  \phi(\lambda) = 0.5 \lambda^2, \ \lambda_0 = \infty; $  is the so-called subgaussian case. \par

  \ We denote the set of all these Young-Orlicz  function as $ \Phi; \ \Phi = \{ \ \phi(\cdot) \}. $ \par

 \ {\bf Definition 1.1.} (See [20], [7].) \par

 \ We say by definition that the {\it centered} random variable (r.v) $ \xi = \xi(\omega) $
belongs to the space $ B(\phi), $ if there exists certain non-negative constant
$ \tau \ge 0 $ such that

$$
\forall \lambda \in (-\lambda_0, \lambda_0) \ \Rightarrow
\max_{\pm} {\bf E} \exp(\pm \lambda \xi) \le \exp[ \phi(\lambda \ \tau) ]. \eqno(1.2)
$$

 \ The minimal non-negative value $ \tau $ satisfying (1.2) for all the values $  \lambda \in (-\lambda_0, \lambda_0), $
is named a $ B(\phi) \ $ norm of the variable $ \xi, $ write

$$
||\xi||B(\phi)  \stackrel{def}{=}
$$

 $$
 \inf \{ \tau, \ \tau > 0: \ \forall \lambda:  \ |\lambda| < \lambda_0 \ \Rightarrow
  \max_{\pm}{\bf E}\exp( \pm \lambda \xi) \le \exp(\phi(\lambda \ \tau)) \}. \eqno(1.3)
 $$

 \vspace{4mm}

 \ These spaces are very convenient for the investigation of the r.v. having a
exponential decreasing tail of distribution, for instance, for investigation of the limit theorem,
the exponential bounds of distribution for sums of random variables,
non-asymptotical properties, problem of continuous and weak compactness of random fields,
study of Central Limit Theorem in the Banach space etc. The detail investigation of these spaces may be found in [2],
 [7]-[8], [9], [12]-[13],  [15], [17], [18], [19], [20] etc. \par

 \ The space $ B(\phi) $ with respect to the norm $ || \cdot ||B(\phi) $ and
ordinary algebraic operations is a rearrangement invariant Banach space in the classical sense [1], chapters 1,2;
 which is in turn isomorphic to the subspace
consisting on all the centered variables of Orlicz's space $ (\Omega, F, {\bf P}), N(\cdot) $
with $ N \ - $ function

$$
N(u) = \exp(\phi^*(u)) - 1, \eqno(1.4)
$$
where
$$
 \phi^*(u) \stackrel{def}{=} \sup_{\lambda \in \Dom(\phi)} (\lambda u - \phi(\lambda)). \eqno(1.5)
$$
 \ The transform $ \phi \to \phi^* $ is called Young-Fenchel, or Legendre transform. The proof of considered
assertion used the properties of saddle-point method and theorem of Fenchel-Moreau:
$$
\phi^{**} = \phi.
$$

  \ Let $  F =  \{  \xi(s) \}, \ s \in S, \ S  $ is an arbitrary set, be the family of somehow
dependent mean zero random variables. The function $  \phi(\cdot) $ may be "constructive" introduced by the formula

$$
\phi(\lambda) = \phi_F(\lambda) \stackrel{def}{=} \max_{\pm} \ln \sup_{s \in S}
 {\bf E} \exp(  \pm \lambda \xi(s)), \eqno(1.6)
$$
 if obviously the family $  F  $ of the centered r.v. $ \{ \xi(s), \ s \in S \} $ satisfies the  so-called
{\it uniform } Cramer's condition:
$$
\exists \mu \in (0, \infty), \ \sup_{s \in S} T_{\xi(s)}(y) \le \exp(-\mu \ y),\ y \ge 0. \eqno(1.7)
$$

 \ Hereafter the symbol $ \ T_{\xi}= T_{\xi}(x) \ $ will be denote the so - called {\it tail function} for the r.v. $ \ \xi: \ $

$$
T_{\xi}(x) \stackrel{def}{=} \max ( {\bf P}(\xi \ge x), \ {\bf P}(\xi \le -x )), \ x \ge 0.
$$

 \ In this case, i.e. in the case the choice the function $ \phi(\cdot) $ by the
formula (1.6), we will call the function $ \phi(\lambda) = \phi_0(\lambda) $
as a {\it natural } function, and correspondingly the function

$$
\lambda \to {\bf E} e^{\lambda \xi}
$$
is named often as a moment generating function  (briefly: MGF)
for the r.v.  $ \xi, $ if of course there exists  in some non-trivial
neighborhood of origin. \par

 \ If for example $  \phi = \phi_2(\lambda) = 0.5 \ \lambda^2, \ \lambda \in R, $ then the r.v. from the space $ B\phi_2 $ are named
subgaussian. This notion was introduced by  J.P.Kahane; see [2]; see also [7] - [8],  [12] - [13], [17]. One can consider
also the case when

$$
\phi_m(\lambda) = m^{-1} \ |\lambda|^m, \ |\lambda| \ge 1, \ m = \const > 1,
$$
as well as a more general case

$$
\phi_{m,L}(\lambda) = m^{-1} \ |\lambda|^m \ L(\lambda), \ |\lambda| \ge 1, \ m = \const > 1,
$$
where $ \ L(\cdot) \ $ is arbitrary positive continuous slowly varying  as $ \ \lambda \to \infty \ $ function,
see [7] - [8]. \par

\vspace{4mm}

 \begin{center}

 \ {\bf B.  Grand Lebesgue \ $ \  G(\psi)  \  $ spaces.} \par

 \end{center}

\vspace{4mm}

  \ Let now  $  \psi = \psi(p), \ p \in [1,b), \ b = \const \in (1, \infty] $ be certain bounded
from below: $ \inf \psi(p) > 0 $ continuous inside the semi - open interval  $ p \in [1,b) $
numerical valued function. We can and will suppose  without loss of generality

 $$
  \inf_{p \in [1,b)} \psi(p) = 1 \eqno(1.8)
$$
and $  \ b = \sup \{ p, \ \psi(p) < \infty  \}, $
 so that   $ \supp \ \psi = [1, b) $ or  $ \supp \ \psi = [1, b]. $ The set of all such a functions will be
denoted by  $ \ \Psi(b) = \{ \psi(\cdot)  \}; \ \Psi := \Psi(\infty).  $ \par

\vspace{4mm}

 {\bf Definition 1.2.} (See [20], [7]-[9].) \par

 \ By definition, the (Banach) Grand Lebesgue Space (GLS)    $  \ G \psi  = G\psi(b) $
consists on all the real (or complex) numerical valued measurable functions
(random variables, r.v.)   $   \  f: \Omega \to R \ $  defined on our probability space and having a finite norm

$$
|| \ f \ || = ||f||G\psi \stackrel{def}{=} \sup_{p \in [1,b)} \left[ \frac{|f|_p}{\psi(p)} \right]. \eqno(1.9)
$$

 \ Here and in what follows the notation $ \  |f|_p  \ $ denotes an ordinary Lebesgue-Riesz $ \ L_p(\Omega) \ $ norm for the r.v. (measurable function) $ \ f: \ $

 $$
 |f|_p  \stackrel{def}{=} \left[ {\bf E} |f|^p   \right]^{1/p}, \ p \ge 1.
 $$

\vspace{4mm}

 \ The function $ \  \psi = \psi(p) \  $ is said to be (also) the {\it  generating function } for this space. \par

 \ If for instance $ \ \psi(p) = \psi_{r}(p) = 1, \ p \in [1,r],  \  $ where $ \ r = \const \in [1,\infty), $ (an extremal case), then the correspondent
 $ \  G\psi_{r}(p)  \  $ space coincides  with the classical Lebesgue - Riesz space

$$
||\xi|| G\psi_{r} = |\xi|_r, \ r \in [1, \infty).
$$

\vspace{4mm}

\  Furthermore,  let now $  \eta = \eta(z), \ z \in S $ be arbitrary family  of random variables  defined on any set $ \ z \in S \ $ such that

$$
\exists b = \const\in (1,\infty], \ \forall p \in [1,b)  \ \Rightarrow  \psi_S(p) := \sup_{z \in S} |\eta(z)|_p  < \infty.
$$
 \ The function $  p \to \psi_S(p)  $ is named as a {\it  natural} function for the  family  of random variables $  S.  $  Obviously,

$$
\sup_{z \in S} ||\eta(z)||G\Psi_S = 1.
$$

 \ The family $ \ S \ $ may consists on the unique r.v., say $  \  \Delta: \ $

$$
\psi_{\Delta}(p):= |\Delta|_p,
$$
if of course  the last function is finite for some value $ \  p = p_0 > 1. \  $\par

\vspace{4mm}

 \begin{center}

 \ {\bf C. Tail inequalities  for both the considered  spaces.} \par

 \end{center}

\vspace{4mm}

 \ {\bf !.} Let $ \ \xi \ $ be non-zero random variable from the space  $ \ B(\phi), \ \phi \in \Phi:  \  K := ||\xi||B(\phi) \ \in (0, \infty). \ $  It follows
immediately by virtue of Tchebychev-Tchernov inequality

$$
T_{\xi}(x) \le  \exp(-\phi^*(x/K)), \ x > 0, \eqno(1.10)
$$
see [7],  [8], [12]. \par

\vspace{4mm}

 \ {\bf 2.} Assume that $ \ \psi(\cdot) \in \Psi(b), \ b = \const \in (1, \infty], \ $ and that the non - zero r.v.  $ \ \eta \ $ belongs to the space
 $ \ G(\psi): \ V  = ||\eta||G\psi \in (0, \infty). \ $ Denote

$$
h(p) = h[\psi](p) \stackrel{def}{=} p \ \ln \psi(p) \eqno(1.11)
$$
 \ Then

$$
T_{\eta}(x) \le  \exp(-h^*(x/V)), \ x > 0. \eqno(1.12)
$$

 \ Both the inequalities (1.10) and (1.12) may be rewritten in the terms of (generalized) Lorentz spaces as follows. Let $ \ S(x), \ x \ge 0 \ $ be
some tail function, i.e. left continuous numerical valued decreasing function $ \  S = S(x) \ $  such that

$$
S(0+) = 1, \ S(\infty) = 0.
$$
 \ A particular cases:

$$
S[\phi](x) :=  \exp(-\phi^*(x)), \ x > 0;
$$

$$
S_{\psi}(x) :=  \exp(-h[\psi]^*(x)), \ x > 0.
$$

 \ Define following, e.g., [1], chapter 4; [5], [10]-[11] the generalized Lorentz quasi - norm $ \  ||\zeta||L(S) \ $
for arbitrary r.v.  $ \ \zeta \ $ as follows.

$$
||\zeta||L(S) \stackrel{def}{=} \sup_{x \ge 0} \left[\ T_{\zeta}(x)/S(x) \ \right]. \eqno(1.13)
$$

 \ We have the following copies of (1.10) and (1.12)

$$
||\xi||L(S[\phi]) \le ||\xi||B(\phi) \eqno(1.14)
$$
and

$$
||\eta||L(S_{\psi}) \le ||\eta||G\psi. \eqno(1.15)
$$

 \ {\sc \ It emerges the following natural  question: to what extent is the converse inequalities  to ones (1.14), (1.15), i.e. when}

$$
||\xi||B(\phi) \le C_1(\phi) \ ||\xi||L(S[\phi]) \eqno(1.16)
$$
{ \sc and}

$$
 ||\eta||G\psi \le  C_2(\psi) ||\eta||L(S_{\psi}) \eqno(1.17)
$$
{\sc for  some finite constants} $ \  C_1(\phi), \ C_2(\psi). \ $ \par

\vspace{5mm}

 {\it Throughout this paper, the letters  $ \ C, C_j(\cdot) \ $ etc. will denote a various positive finite
 constants which may differ from one formula to the next even within a single string
 of estimates and which does not depend on the essentially variables  $  \ p, x, \lambda, y, u \ $ etc. \par

 \ We make no attempt to obtain the best values for these constants.}\par

\vspace{5mm}

 \ We represent  for beginning some results concerning both the inverse estimates  (1.6), (1.17). \par

 \  It is proved in particular in the article [8], Theorem 4.1, that if

$$
 R[\psi] := \sup_{p \in [1, \infty)} \left[ \ h^{* `} [\psi](p) \   \right]^{1/p} < \infty, \eqno(1.18)
$$
then

$$
 ||\eta ||G\psi \le 2 \   R[\psi] \  e^{1/e} \ ||\eta||L(S_{\psi}) < \infty, \eqno(1.19)
$$
i.e. the inequality  (1.17). \par

\vspace{4mm}

 \ Let us turn our attention on the relation (1.16). We show here briefly the result from the aforementioned article [8]. \par

 \ In order to carry out this,  we assume here $ \ \lambda_0 = \lambda_0[\phi] = \infty, \ $  and define for any positive finite constant value $ \ C_1 \ $
the function

$$
\theta[\phi](\lambda) =
\theta(\lambda) \stackrel{def}{=} \frac{C_1}{\lambda \ \phi^{*'}(\lambda)} \eqno(1.20)
$$
for all the sufficiently greatest values $  \lambda: \ \lambda \ge e, \ $ (say).  \ We introduce also the following integral

$$
Z[\phi](\lambda) = Z(\lambda) := \int_0^{\infty} e^{ - \theta(\lambda) \ \phi^*(x) } \ dx. \eqno(1.21)
$$

\vspace{3mm}

\ Suppose

$$
\exists C_3 = C_3(C_1) = \const < \infty, \ \ \forall \lambda > e \ \Rightarrow Z[\phi](\lambda) \le  \exp \phi^*(C_3 \lambda). \eqno(1.22)
$$

 \ Then the inequality (1.16) holds true. \par

\vspace{3mm}

 {\bf Examples.} \par

\vspace{3mm}

{\it  Example 1.}  Let $ \ L = L(\lambda), \ x > 0 \ $ be positive twice continuous differentiable slowly varying at infinity regular in the following sense

$$
\lim_{\lambda \to \infty}  \frac{L(\lambda/L(\lambda))}{L(\lambda)} = 1
$$
function. Define also for sufficiently greatest values $ \ \lambda, \ $  say for $ |\lambda| \ge 1, \ $  the  function of the form

$$
\phi_{m,L}(\lambda) \stackrel{def}{=} m^{-1} \ |\lambda|^m \ L^{1/q}(|\lambda|^m), \
$$
$ \ m = \const > 1, \ q = m/(m-1); $ and as usually

$$
\phi_{m,L}(\lambda) \stackrel{def}{=} C(m,L) \ \lambda^2, |\lambda| < 1.
$$

 \ It is proved in [8] in particular that for these $ \ \Phi \ - $ function  the equality (1.16) there holds.
Namely, the inclusion

$$
\xi \in B(\phi_{m,L}) < \infty, \ \xi \ne 0
$$
is quite equivalent to the following tail estimate

$$
\exists C_2 \in (0, \infty) \ \Rightarrow T_{\xi}(x) \le \exp \left( - C_2 \ q^{-1} \ x^q \  L^{-(q-1)} \left(x^{q-1} \right)  \right), \  x \ge 1.
$$

\vspace{3mm}

{\it  Example 2.} Define the other Grand Lebesgue Space space $ \ G \psi_{C,\beta}   \ $ of random variables with
correspondent generating function

$$
\psi_{C,\beta}(p) \stackrel{def}{=} e^{C \ p^{\beta}}, \ p \in [1,\infty),
$$
where $ \  \beta,C = \const > 0. \ $ \ It is proved in [8] in particular that for these $ \ \Psi \ - $ function  the equality (1.17) again there holds. \par

 \ In detail, the following implication there holds

$$
\xi \in \cup_{C >0} G\psi_{C,\beta} \Leftrightarrow
$$

$$
 \exists K \in (0,\infty), \ T_{\xi}(x) \le \exp\left(  - K \ \left(\ln (x+1)^{1 +1/ \beta} \right) \right), \ x \ge 0.
$$

\vspace{4mm}

 \ Note that in general case the MGF for arbitrary r.v. $ \  \nu \ $ with

$$
||\eta||G \psi_{C,\beta} \in (0, \infty)
$$
does not exists; on the other words this variable does not satisfy the Cramer's condition. \par

\vspace{3mm}

{\it  Example 3.} Define the following $ \ \Psi \ $ function

$$
\psi_m(p) \stackrel{def}{=} p^{1/m}, \ p \in [1,\infty); \ m = \const > 0.
$$
 \ The non - zero r.v. $ \  \xi \ $ belongs to the space $ \  G\psi_m:  \ $

$$
||\xi||G \psi_m = \sup_{p \ge 1} \left[ \ \frac{|\xi|_p}{p^{1/m}} \ \right] < \infty
$$
  if and only if

$$
\exists C_3 \in (0,\infty) \ \Rightarrow T_{\xi}(x) \le \exp \left( - C_3 \ x^m \right), \ x \ge 1.
$$

\vspace{4mm}

\ {\bf Our purpose  in this short article is to show  by means of building of suitable counterexamples that both the estimates (1.16) and (1.17)
 are not true if the conditions correspondingly relations (1.18) and (1.21) are not satisfied.  } \par

 \ Obtained here results generalized ones in [8]. \par

\vspace{4mm}

\section{Main results: counterexamples.}

 \vspace{4mm}

 \ {\bf A. Grand Lebesgue Space case. } \par

 \vspace{4mm}

\ {\bf I.}  \ Let us introduce the following $ \ \Psi \ $  function

 $$
 \psi^{b, \gamma}(p) = (b-p)^{-\gamma/b},  \ p \in [1,b), \eqno(2.1)
 $$
where  $ b = \const > 1, \ \gamma = \const > 0.$ \par

 \ Suppose the non-zero r.v. $ \ \xi \ $ belongs to the space $ \ G\psi^{b, \gamma}; \ $  one can conclude without loss of
 generality $ \  ||\xi|| \ G\psi^{b, \gamma} = 1. \  $ We find my means of simple calculations based on the  estimate  (1.12)

$$
T_{\xi}(x) \le C(b,\gamma) \ x^{-b} \ [\ln x]^{\gamma}, \ x \ge e. \eqno(2.2)
$$

\vspace{4mm}

 \ {\bf II.}  Conversely, consider the r.v. $ \ \zeta \ $ with tail behavior of a type (2.2):

$$
T_{\zeta}(x) =  x^{-b} \ [\ln x]^{\gamma}, \ x \ge x_0 = \const \ge e. \eqno(2.3)
$$
 \ We have as $ \ p \to b - 0 $

$$
C_3^{-1}  \ p^{-1} \ |\zeta|_p^p  \sim  \int_1^{\infty} x^{p-1} \ x^{-b} \ (\ln x)^{\gamma} \ dx =
$$

$$
\int_0^{\infty} e^{ -(b-p) \ y } \ y^{\gamma} \ dy  = \frac{\Gamma(\gamma + 1)}{(b-p)^{\gamma + 1}};
$$

$$
|\zeta|_p \asymp (b - p)^{-(\gamma + 1)/b}, \ p \in [1,b). \ \eqno(2.4)
$$
 \ Thus, the r.v. $ \ \zeta \ $  belongs to the space $ \   G\psi^{b, \gamma + 1},  \ $ but it
 does not belongs to the space $ \   G\psi^{b, \gamma}:  \ $

$$
||\zeta|| G\psi^{b, \gamma} = \infty. \eqno(2.5)
$$

\vspace{4mm}

 \ {\it For illustration:}  consider the space $ \ G\psi_{(r)}; \ $  if for some r.v. $ \ \zeta \ ||\zeta||  G\psi_{(r)} = 1, \ r = \const \ge 1, \ $ or equally
 $ \ |\zeta|_r = 1, \ $ then

 $$
 T_{\zeta}(x) \le \frac{1}{x^r}, \ x \ge 1;
 $$
but the inverse conclusion is obviously not true. \par

\vspace{4mm}

 \ {\bf B. $ \ B(\phi) \ $ space case. } \par

 \vspace{4mm}

\ {\bf I.}  \ Let us introduce the following $ \ \Phi \ $  function

$$
\phi_{b, \gamma}(\lambda) := \gamma \ \ln \left[ \frac{b}{b-|\lambda|} \right]; \ |\lambda| \in [1,b),  \eqno(2.6)
$$
where  $ \ \gamma = \const > 0; \ b = \const > 1.  \ $ \par
 \ If the non-zero random variable $ \ \tau \ $ belongs to the space $ \  B( \phi_{b, \gamma}), \  $ for instance
when

$$
||\tau|| B(\phi_{b, \gamma}) = 1,
$$
then

$$
T_{\tau}(x)  \le C \ x^{\gamma} \ e^{ - b \ x  }, \ x \ge e.  \eqno(2.7)
$$

\vspace{4mm}

\ {\bf II.} Inversely, let the r.v. \ $ \theta \ $  has a  tail behavior of a type (2.7):

$$
T_{\theta}(x) = C \ x^{\gamma} \ e^{ - b \ x  }, \ x \ge e.  \eqno(2.8)
$$

 \ Then  as $ \ \lambda \to  b - 0,  \ \lambda \in [1,b)  $

$$
{\bf E} e^{\lambda \ \theta}  \sim \lambda \int_0^{\infty} e^{ \ - (b - \lambda) x \ }  \ x^{\gamma } \ dx =
$$

$$
\frac{\lambda \ \Gamma(\gamma +1)}{(b - \lambda)^{\gamma + 1}} \asymp \left( \ \frac{b}{b - \lambda} \  \right)^{\gamma + 1}. \eqno(2.9)
$$

 \  Therefore, the r.v. $ \ \theta \ $ does not belongs to the space  $ \  B(\phi_{b, \gamma}). \ $ \par

\vspace{4mm}

\section{Concluding remarks.}

\vspace{4mm}

{\bf A.} \ It is interest by our opinion to obtain the generalization of  results of this report into the multidimensional
case, i.e. into random vectors, alike in the articles [8], [13]. \par

\vspace{4mm}

{\bf B. }  We mention  even briefly an important  possible application of obtained results: a Central Limit Theorem in Banach
spaces, in the spirit of [15], [14], [16], [12], section 4.1. \par

\vspace{4mm}

{\bf C.}  Authors have not  the correspondent counterexamples in the case when $ \ b = \infty. \ $ This is an open problem. \par

\begin{center}

 \vspace{6mm}

 {\bf References.}

\end{center}
 \vspace{4mm}

{\bf 1. Bennet C., Sharpley R.}  {\it  Interpolation of operators.} Orlando, Academic
Press Inc., (1988). \\

 \vspace{3mm}

{\bf 2.  Buldygin V.V., Kozachenko Yu.V. }  {\it Metric Characterization of Random
Variables and Random Processes.} 1998, Translations of Mathematics Monograph, AMS, v.188. \\

 \vspace{3mm}

 {\bf 3. A. Fiorenza.}   {\it Duality and reflexivity in grand Lebesgue spaces. } Collect. Math.
{\bf 51,}  (2000), 131-148. \\

 \vspace{3mm}

{\bf  4. A. Fiorenza and G.E. Karadzhov.} {\it Grand and small Lebesgue spaces and
their analogs.} Consiglio Nationale Delle Ricerche, Instituto per le Applicazioni
del Calcoto Mauro Picone, Sezione di Napoli, Rapporto tecnico 272/03, (2005).\\

 \vspace{3mm}

{\bf 5. Grafakos, Loukas.} {\it Classical Fourier analysis.} Graduate Texts in Mathematics, 249, (2nd ed.),
Berlin, New York: Springer-Verlag, (2008), ISBN 978-0-387-09431-1, MR 2445437, doi:10.1007/978-0-387-09432-8. \\

\vspace{3mm}

{\bf 6.  T. Iwaniec and C. Sbordone.} {\it On the integrability of the Jacobian under minimal
hypotheses. } Arch. Rat.Mech. Anal., 119, (1992), 129-143. \\

 \vspace{3mm}

{\bf 7. Kozachenko Yu. V., Ostrovsky E.I. }  (1985). {\it The Banach Spaces of random Variables of subgaussian Type. } Theory of Probab.
and Math. Stat. (in Russian). Kiev, KSU, 32, 43-57. \\

\vspace{3mm}

{\bf 8. Kozachenko Yu.V., Ostrovsky E., Sirota L.}  {\it Relations between exponential tails, moments and
moment generating functions for random variables and vectors.} \\
arXiv:1701.01901v1 [math.FA] 8 Jan 2017 \\

\vspace{3mm}

{\bf 9.  E Liflyand, E Ostrovsky, L Sirota.} {\it Structural properties of bilateral grand Lebesgue spaces.}
Turkish Journal of Mathematics, (2010), {\bf 34, } \ (2), 207-220. \\

\vspace{3mm}

{\bf 10. G. Lorentz.} {\it Some new function spaces.}  Annals of Mathematics, {\bf 51}, (1950), pp. 37-55.\\

\vspace{3mm}

{\bf 11. G. Lorentz.} {\it On the theory of spaces } $ \ \Lambda. \ $  Pacific Journal of Mathematics, {\bf 1,}  (1951), pp. 411-429.\\

\vspace{3mm}

{\bf 12. Ostrovsky E.I. } (1999). {\it Exponential estimations for Random Fields and its
applications,} (in Russian). Moscow-Obninsk, OINPE. \\

 \vspace{3mm}

{\bf 13. Ostrovsky E. and Sirota L.} {\it Vector rearrangement invariant Banach spaces
of random variables with exponential decreasing tails of distributions.} \\
 arXiv:1510.04182v1 [math.PR] 14 Oct 2015 \\

 \vspace{3mm}

{\bf 14. Ostrovsky E. and Sirota L.}  {\it Non-asymptotical sharp exponential estimates
for maximum distribution of discontinuous random fields. } \\
 arXiv:1510.08945v1 [math.PR] 30 Oct 2015 \\

 \vspace{3mm}

 {\bf 15. Ostrovsky E., Rogover E. } {\it Exact exponential bounds for the random field
 maximum distribution via the majorizing measures (generic chaining).} \\
 arXiv:0802.0349v1 [math.PR] 4 Feb 2008 \\

 \vspace{3mm}

{\bf 16. Ostrovsky E. and Sirota L. } {\it   Entropy and Grand Lebesgue Spaces approach for the problem  of
Prokhorov-Skorokhod continuity of discontinuous random fields. }\\
arXiv:1512.01909v1 [math.Pr] 7 Dec 2015 \\

\vspace{3mm}

{\bf 17. Ostrovsky E. and Sirota L. } {\it Moment Banach spaces: Theory and applications.}
2007, Holon Institute of Technology,
HAIT Journal of Science and Engineering C, Volume 4, Issues 1-2, pp. 233-262.\\

 \vspace{3mm}

 {\bf 18. Ostrovsky E. and Sirota L. } {\it  Fundamental function for Grand Lebesgue Spaces.  } \\
 arXiv:1509.03644v1  [math.FA]  11 Sep 2015 \\

 \vspace{3mm}

 {\bf 19. Ostrovsky E.} {\it Generalization of the norms of Buldygin-Kozachenko and  Central Limit Theorem in Banach space.}
Probab. Theory Applications, 1982, V. 27  Issue 3, p. 618, (in Russian). \\

\vspace{3mm}

{\bf 20.  O.I.Vasalik, Yu.V.Kozachenko, R.E.Yamnenko.} {\it $ \ \phi \ - $ subgaussian  random processes. } Monograph, Kiev, KSU,
2008;  (in Ukrainian). \\

\vspace{3mm}

\end{document}